\documentclass[11pt]{article}
\usepackage{a4}
\usepackage{color,graphicx,bbm}

\usepackage{amssymb,amsmath,amsthm}

\newcommand{\Z}{{\bf Z}}

\DeclareMathOperator{\st}{st}
\DeclareMathOperator{\lk}{lk}
\DeclareMathOperator{\antist}{ast}
\DeclareMathOperator{\dd}{dist}
\DeclareMathOperator{\g}{g}
\DeclareMathOperator{\f}{f}
\DeclareMathOperator{\Rig}{Rig}
\DeclareMathOperator{\Ker}{Ker}
\DeclareMathOperator{\rank}{rank}
\DeclareMathOperator{\Imm}{Im}

\makeatletter
\def\imod#1{\allowbreak\mkern10mu({\operator@font mod}\,\,#1)}
\makeatother

\theoremstyle{plain}
\newtheorem{theorem}{Theorem}[section]
\newtheorem{proposition}[theorem]{Proposition}
\newtheorem{lemma}[theorem]{Lemma}
\newtheorem*{lemma*}{Auxiliary Lemma}
\newtheorem{observation}[theorem]{Observation}
\newtheorem{corollary}[theorem]{Corollary}

\newtheorem{problem}[theorem]{Problem}
\newtheorem{thmRemark}[theorem]{Remark}

\theoremstyle{definition}

\newtheorem{remark}[theorem]{Remark}
\newtheorem{remarks}[theorem]{Remarks}
\newtheorem{example}[theorem]{Example}
\newtheorem*{acknowledgement}{Acknowledgement}

\title{A characterization of simplicial polytopes with $\g_2=1$}
\author{Eran Nevo\footnote{Department of Mathematics, Cornell
University, Ithaca USA, E-mail address:
eranevo@math.cornell.edu. Research partially supported by an NSF grant DMS-0757828.} \ and Eyal Novinsky\footnote{Institute of Mathematics, The Hebrew
University, Jerusalem Israel, E-mail address:
eyalnov@math.huji.ac.il}}

\begin{document}
\maketitle
\begin{abstract}
Kalai proved that the simplicial polytopes with $\g_2=0$ are the stacked polytopes. We characterize the $\g_2=1$ case.

Specifically, we prove that every simplicial $d$-polytope ($d\geq 4$) which is prime and with $\g_2=1$ is combinatorially equivalent either to a free sum of two simplices whose dimensions add up to $d$ (each of dimension at least $2$), or to a free sum of a polygon with a $(d-2)$-simplex. Thus, every simplicial $d$-polytope ($d\geq 4$) with $\g_2=1$ is combinatorially equivalent to a polytope obtained by stacking over a polytope as above. Moreover, the above characterization holds for any homology $(d-1)$-sphere ($d\geq 4$) with $\g_2=1$, and our proof takes advantage of working with this larger class of complexes.
\end{abstract}

\section{Introduction and results}\label{Introduction}
Let $\f_i(K)$ denote the number of $i$-dimensional faces in a simplicial complex $K$. In particular, $\f_0$ counts vertices and $\f_1$ counts edges. Let $$\g_2(K):= \f_1(K)-d \f_0(K)+\binom{d+1}{2}$$ where $d$ is the maximal size of a face of $K$, i.e. $d$ equals the dimension of $K$ plus $1$.
(This notation is standard in face-vector theory, see e.g. \cite{Kalai:Aspects-94} for details.)
The well known Lower Bound Theorem (LBT) proved by Barnette \cite{Barnette-LBT,Barnete-LBTfacets,BarnetteManifolds}, asserts that if $K$ is the boundary complex of a simplicial $d$-polytope, or more generally a finite triangulation of a connected compact $(d-1)$-manifold without boundary, where $d\geq 3$, then $\g_2(K)\geq 0$.
Kalai considered several generalizations of this result, including to homology manifolds, and characterized the case of equality \cite{Kalai-LBT}. To state his result, define stacked polytopes:
a \emph{stacking} is the operation of adding a pyramid over a facet of a given simplicial polytope.
A polytope is \emph{stacked} if it can be obtained from a simplex by repeating the stacking operation (finitely many times, may be zero).
We will make use of the following result:

\begin{theorem}\label{thm:stacked,Kalai}\cite{Barnette-LBT,BarnetteManifolds} and \cite[Theorems 6.2 and 7.1]{Kalai-LBT}
Let $d\geq 4$, and let $K$ be the boundary complex of a simplicial $d$-polytope, or more generally a homology $(d-1)$-manifold. Then $\g_2(K)\geq 0$ and equality holds iff $K$ is combinatorially isomorphic to the boundary complex of a stacked $d$-polytope.
\end{theorem}
Kalai's proof is based on results from rigidity theory, to be discussed later; see also Gromov \cite{Gromov} for a proof of the nonnegativity of $\g_2$.

A subset $F$ of the vertices of a simplicial complex $K$ is called a \emph{missing face of $K$} if $F\notin K$ and all proper subsets of $F$ are in $K$; if $F$ has size $d$ it is called a \emph{missing (d-1)-face}, and if in addition $K$ is $(d-1)$-dimensional then we say that $F$ is a \emph{missing facet} of $K$.
A simplicial $d$-polytope is called \emph{prime} if its boundary complex contains no missing $(d-1)$-faces and is not the boundary complex of a simplex. Similarly, prime homology spheres are defined as homology spheres with no missing facets. 
We call a simplicial prime polytope a \emph{prime polytope} for short. For example, for the two $3$-polytopes bipyramid and octahedron, the first is not prime as is has a missing $2$-face, while the second is prime, as it has no missing $2$-face.

In \cite[Theorem 3.10]{Kalai:Aspects-94} Kalai claimed that there exists a function $u(d,b)$ such that if the boundary complex $K$ of a prime $d$-polytope ($d\geq 4$) satisfies $\g_2(K)=b$ then $\f_0(K)\leq u(d,b)$. We provide a counterexample (Example \ref{ex:big g1}). First, let us fix some notation:
the boundary complex of a simplicial polytope $P$ is denoted by $\partial (P)$, or simply by $\partial P$.
The \emph{free sum} of two polytopes $P,Q$, denoted by $P\star Q$, is defined as the convex hull of their union when $P$ and $Q$ are embedded in orthogonal spaces with the origin in the interior of both. Indeed, the combinatorial type of $P\star Q$ is well defined: its boundary complex is the join $\partial P*\partial Q$, where the \emph{join} of two simplicial complexes $K,L$ is the collection of disjoint unions $\{A\uplus B ~:~A\in K,\ B\in L\}$.
A direct computation shows:
\begin{example}\label{ex:big g1}
Let $C_n$ be a $2$-polytope with $n$ vertices and let $\sigma^m$ be the $m$-dimensional simplex.  Then for every $d\geq 4$ and any $n\geq 3$, $C_n\star \sigma^{d-2}$ is a prime $d$-polytope with  $\g_2(\partial(C_n \star \sigma^{d-2}))=1$.
\end{example}

Our main result characterizes the prime polytopes with $\g_2=1$:

\begin{theorem}\label{thm:Main}
Let $d\geq 4$, and let $K$ be the boundary complex of a prime $d$-polytope, or more generally a prime homology $(d-1)$-sphere. Assume that $\g_2(K)=1$. Then $K$ is combinatorially isomorphic to either the join of boundary complexes of two simplices whose dimensions add up to $d$ (each simplex is of dimension at least $2$), or the join of the boundary complexes of a convex polygon and a $(d-2)$-simplex.
\end{theorem}
Note that any simplicial polytope can be (uniquely) presented as a connected sum of prime polytopes and simplices, and similarly for homology spheres, and that $\g_2$ of a connected sum is the sum of $\g_2$'s of its components.
(Recall that the \emph{connected sum} of two disjoint simplicial complexes of equal dimension is the operation of identifying by a bijection the vertices in a facet of one with the vertices in a facet of the other, identifying the faces they form accordingly, and later deleting the identified facet. Thus, the connected sum of homology spheres is a homology sphere, by an easy Mayer-Vietoris argument and Alexander duality. For polytopes, after suitable projective transformations of each, which of course preserve their combinatorial structure, the connected sum, which is gluing along a facet of each, can be made convex too.)
Thus, by Theorems \ref{thm:stacked,Kalai} and \ref{thm:Main} we conclude that:

\begin{corollary}
Let $d\geq 4$, and let $K$ be the boundary complex of a $d$-polytope, or a homology $(d-1)$-sphere, with $\g_2(K)=1$. Then $K$ is combinatorially isomorphic to the boundary complex of a polytope obtained by repeated stacking, starting from either the free sum of two simplices whose dimensions add up to $d$ (each simplex is of dimension at least $2$), or from the free sum of a polygon and a $(d-2)$-simplex.
\end{corollary}
This result can be compared with Perles' characterization of polytopes $P$ with $\g_1(P):=\f_0(P)-(\dim(P)+1)$ at most $2$ \cite[Chapter 6]{Grunbaum:ConvexPolytopes-03} and with Mani's result that triangulated spheres with $\g_1\leq 2$ are polytopal \cite{Mani}.
We do not know of a characterization of simplicial polytopes with $\g_2=2$.

The proof of Theorem \ref{thm:Main} is based on rigidity theory for graphs, introduced in \cite{Asi-Roth1,Asi-Roth2}.
Working with homology spheres, rather than with simplicial polytopes, greatly simplifies the proof; in particular see the proof of Proposition \ref{prop:d=4,no missing trianle}.

This paper is organized as follows: in Section \ref{sec:Background} we give the necessary background for polytopes and homology spheres, and develop the needed results in rigidity theory of graphs. In Section \ref{sec:Proofs} we prove Theorem \ref{thm:Main} and discuss some extensions of it and related open problems.

\section{Background}\label{sec:Background}
\paragraph{Polytopes and homology spheres.}
For unexplained terminology we refer to textbooks on polytopes, e.g. \cite{Grunbaum:ConvexPolytopes-03,Ziegler}, and on simplicial homology, e.g. \cite{Munkres}.

In this paper a \emph{simplicial complex} means a finite abstract simplicial complex, i.e. a finite collection of finite sets closed under inclusion. 
The \emph{$i$-skeleton} of a simplicial complex $K$ is $K_{\leq i}=\{F\in K: |F|\leq i+1\}$. The \emph{graph} of $K$ is $K_{\leq 1}$. Let $K_i:=\{F\in K: |F|=i+1\}$. 
The \emph{link} of a face $F$ in a $K$ is $\lk(F)=\lk(F,K)=\{T\in K: \ T\cap F=\emptyset, \ T\cup F\in K\}$, its \emph{closed star} is $\overline{\st(F)}=\overline{\st(F,K)}=\{T\in K: \ T\cup F\in K\}$,
its \emph{antistar} is $\antist(F)=\antist(F,K)=\{T\in K: \ T\cap F=\emptyset\}$; they are simplicial complexes as well. The (open) \emph{star} of $F$ is the collection of sets $\st(F)=\st(F,K)=\{T\in K: \ F\subseteq T\}$.

Note that for a vertex $v$ in a simplicial polytope $P$, its vertex figure $P/v$ satisfies $\partial(P/v)=\lk(v,\partial P)$.

A \emph{homology sphere} is a simplicial complex $K$ such that for every face $F$ in $K$ (including the empty set), and for every $0\leq i$ there is an isomorphism of reduced homology groups $\tilde{H}_i(\lk(F,K);\Z)\cong \tilde{H}_i(S^{\dim(K)-|F|};\Z)$ where $S^m$ denotes the $m$-dimensional sphere and $\Z$ the integers (actually any fixed coefficients ring works for Theorem \ref{thm:Main}). In particular, a boundary complex of a simplicial polytope is a homology sphere; however there are many non-polytopal examples of homology spheres, e.g. \cite{Kalai-manyspheres}.
Alexander Duality holds for homology spheres, e.g. \cite[Chapter 8, \S 71]{Munkres}, cited below. Denote by $||K||$ a geometric realization of
a simplicial complex $K$, and for a subcomplex $A$ of $K$ let $||A||$ denote the subset of $||K||$ induced by the inclusion $A\subseteq K$. Let $\tilde{H}^k$ denote reduced $k$'th  cohomology (say with integer coefficients).
\begin{theorem}\label{thm:Alexander Duality}(Alexander Duality)
Let $A$ be a proper nonempty subcomplex of a homology $n$-sphere $K$. Then for every $k$, $\tilde{H}^k(A)\cong \tilde{H}_{n-k-1}(||K||-||A||)$.  

In particular (we will use only these facts in the sequel),
such $A$ is never a homology $n$-sphere, and if $A$ is a homology $(n-1)$-sphere then $K-A$ has two connected components, and $A$ is their common boundary.
\end{theorem}
A \emph{homology ball} is an acyclic simplicial complex $K$ where the link of every face $F$ in $K$ is either acyclic or has the homology of a sphere of dimension $\dim{K}-|F|$, and the faces with acyclic links form a homology $(\dim{K}-1)$-sphere, called the \emph{boundary} of $K$. Note that if $K$ is a homology sphere and $v$ a vertex in $K$ then $\antist(v,K)$ is a homology ball of the same dimension as $K$.

We will use the following known fact \cite[Corollary 70.3]{Munkres}:
\begin{proposition}\label{prop:facet-connected}
Let $K$ be a homology ball (or more generally a connected homology manifold with or without boundary) of dimension $\geq 1$. Then $K$ is facet connected, i.e. for any two maximal faces of $S,T$ of $K$ there is a sequence of facets $(S=F_0,F_1,...,F_i=T)$ in $K$ such that $F_j\cap F_{j-1}$ is a codimension $1$ face of $K$ for every $1\leq j\leq i$.
\end{proposition}
Note that for $K=\partial P$, $P$ a simplicial polytope, and $v\in K_0$, one can show that $\antist(v,K)$ is facet connected by a line shelling through $v$.

\begin{proposition}\label{prop:edge-contraction-homology-sphere}
Let $K$ be a homology $d$-sphere ($d\geq 1$), and $\{u,v\}\in K_1$ an edge satisfying the \emph{link condition} $\lk(\{u,v\},K)=\lk(u,K)\cap \lk(v,K)$.
Let $K':=(K-\st(u))\cup (\{v\}*\antist(v,\lk(u)))$.
Then $K'$ is a homology $d$-sphere.
\end{proposition}
\begin{proof}
The proof uses an easy Mayer-Vietoris argument. Let $T\in K'$.

Case 1: $T\in K-\st(u)$, $T\notin \lk(u)$. Then $\lk(T,K')=\lk(T,K-\st(u))=\lk(T,K)$ where the last equality follows from $T\notin \partial(K-\st(u))=\lk(u)$. Thus $\tilde{H}_i(\lk(T,K');\Z)\cong \tilde{H}_i(S^{\dim(K')-|T|};\Z)$ for all $i$.

Case 2: $T\in \{v\}*\antist(v,\lk(u))$, $T\notin \lk(u)$. Then $v\in T$ and $\lk(T,K')=\lk(T\setminus \{v\}, \antist(v,\lk(u))= \lk(T\setminus \{v\}, \lk(u))$ where the last equality follows from $T\setminus \{v\}\notin \partial(\antist(v,\lk(u)))=\lk(\{u,v\})$. As $\dim(K')=\dim(\lk(u))+1$ and $\lk(u)$ is a homology sphere, one gets $\tilde{H}_i(\lk(T,K');\Z)\cong \tilde{H}_i(S^{\dim(K')-|T|};\Z)$ for all $i$.

Case 3: $T\in \lk(u)$. Then $\lk(T,K')=\lk(T,K-\st(u))\cup_{\lk(T,\lk(u))} \lk(T,\{v\}*\antist(v,\lk(u)))$. Note that both $K-\st(u)=\antist(u,K)$ and $\{v\}*\antist(v,\lk(u)))$ are homology balls with $T$ contained in their (common) boundary. Thus, plugging $\tilde{H}_i(\lk(T,K-\st(u))=\tilde{H}_i(\lk(T,\{v\}*\antist(v,\lk(u))))=0$ for every $i$ into the Mayer-Vietoris long exact sequence for $\lk(T,K')=\lk(T,K-\st(u))\cup_{\lk(T,\lk(u))} \lk(T,\{v\}*\antist(v,\lk(u)))$ yields $\tilde{H}_i(\lk(T,K')\cong\tilde{H}_{i-1}(\lk(T,\lk(u)))$. Hence, for every $i$, $\tilde{H}_i(\lk(T,K');\Z)\cong \tilde{H}_i(S^{\dim(K')-|T|};\Z)$.

Thus, by definition, $K'$ is a homology $d$-sphere.
\end{proof}
\paragraph{Rigidity.}
The presentation here is based mainly on Kalai's \cite{Kalai-LBT}.
Let $G=(V,E)$ be a simple graph. Let $\dd(a,b)$ denote Euclidian distance
between points $a$ and $b$ in Euclidian space.
A \emph{$d$-embedding} is a map $f:V\rightarrow \mathbb{R}^{d}$. It is called \emph{rigid}
if there exists an $\varepsilon>0$ such that if $g:V\rightarrow
\mathbb{R}^{d}$ satisfies $\dd(f(v),g(v))<\varepsilon$ for every
$v\in V$ and $\dd(g(u),g(w))=\dd(f(u),f(w))$ for every $\{u,w\}\in E$,
then $\dd(g(u),g(w))=\dd(f(u),f(w))$ for every $u,w\in V$.
Loosely speaking, $f$ is rigid if any perturbation of it which
preserves the lengths of the edges is induced by an isometry of $\mathbb{R}^d$.
$G$ is called \emph{generically d-rigid} if the set of its rigid
$d$-embeddings is open and dense in the topological vector space
of all of its $d$-embeddings.
Given a $d$-embedding $f:V\rightarrow \mathbb{R}^{d}$, a \emph{stress}
w.r.t. $f$ is a function $w:E\rightarrow \mathbb{R}$ such that for
every vertex $v\in V$
$$\sum_{u:\{v,u\}\in E}w(\{v,u\})(f(v)-f(u)) =0.$$
$G$ is called \emph{generically d-stress free} if the set of its
$d$-embeddings which have a unique stress ($w=0$) is open and dense
in the space of all of its $d$-embeddings.

Rigidity and stress freeness can be related as follows: Let
$V=[n]$, and let $\Rig(G,f)$ be the $dn\times |E|$ matrix
associated with a $d$-embedding $f$ of $V(G)$ defined as follows:
for its column corresponding to $\{v<u\}\in E$ put the vector
$f(v)-f(u)$ (resp. $f(u)-f(v)$) at the entries of the rows
corresponding to $v$ (resp. $u$) and zero otherwise. $G$ is
generically $d$-stress free iff the kernel $\Ker(\Rig(G,f))=\{0\}$ for a generic $f$
(i.e. for an open dense set of embeddings). $G$ is generically
$d$-rigid iff the images $\Imm(\Rig(G,f))=\Imm(\Rig(K_{V},f)$ for a generic $f$,
where $K_{V}$ is the complete graph on $V=V(G)$. The dimensions of
the kernel and image of $\Rig(G,f)$ are independent of the generic
$f$ we choose; $\Rig(G,f)$ is the \emph{rigidity matrix} of
$G$, denoted by $\Rig(G,d)$ for a generic $f$.
For the complete graph, one computes $\rank(\Rig(K_{V},d))=dn-\binom{d+1}{2}$ (see
Asimov and Roth \cite{Asi-Roth1} for more details).
In particular, if $G$ is generically $d$-rigid then $\g_2(G)$ is the dimension of $\Ker(\Rig(G,d))$.
We say that an edge $\{u,v\}$ \emph{participates} in a stress $w$ if $w(\{u,v\})\neq 0$, and that a vertex $v$ \emph{participates} in $w$ if there exists a vertex $u$ such that the edge $\{u,v\}$ participates in $w$.
A \emph{generic d-stress} is a stress w.r.t. a generic $d$-embedding.

We need the following known results for the proof of Theorem \ref{thm:Main}:
\begin{lemma}\label{lemma:Cone}(Cone Lemma \cite[Theorem 5]{Whiteley-cones}, also \cite[Theorem 1.3, Corollary 1.2]{TWW2})
Let $C(G)$ be the graph of the cone over a graph $G$, i.e. $C(G)=(\{u\}*G)_{\leq 1}$ where $u\notin G$. Then for every $d>0$, 

(1) $G$ is generically $d$-rigid iff $C(G)$ is generically $(d+1)$-rigid.

(2) $\Ker(\Rig(C(G),d+1))\cong\Ker(\Rig(G,d))$ as real vector spaces. Moreover, $u$ participates in a generic stress of $C(G)$, provided that $\Ker(\Rig(G,d))\neq \{0\}$.
\end{lemma}
\begin{remark}
The `moreover part' does not appear explicitly in \cite{Whiteley-cones,TWW2} but is clear for generic embeddings from the isomorphism constructed there.
\end{remark}

\begin{lemma}\label{lemma:Gluing}(Gluing Lemma \cite{Asi-Roth2})
Let $G_{i}=(V_{i},E_{i})$ be generically $d$-rigid graphs, $i=1,2$, such that
$G_{1}\cap G_{2}$ has at least $d$ vertices. Then $G_{1}\cup G_{2}$ is
generically $d$-rigid.
\end{lemma}

The following is well known, e.g. by Cauchy's rigidity theorem for polytopes, or by Gluck \cite{Gluck} for triangulated $2$-spheres:
\begin{lemma}\label{thm:rigid stress free 2-spheres}
Let $G$ be the graph of a convex $3$-polytope (or of a homology $2$-sphere. These two families of graphs coincide). Then $G$ is generically $3$-rigid and $3$-stress free.
\end{lemma}
Using this fact as the base of induction, the Cone Lemma and essentially the Gluing Lemma for the induction step, Kalai \cite{Kalai-LBT} proved:
\begin{lemma}\label{lem:kalai rigid spheres}
Graphs of homology $(d-1)$-spheres are generically $d$-rigid for $d\geq 3$.
\end{lemma}

\begin{lemma}\label{lemma:kalai stacked links}(\cite[Theorems 7.1 and 9.3]{Kalai-LBT})
Let $d>4$, and let $K$ be the boundary complex of a $d$-polytope, or a homology $(d-1)$-sphere. If for every vertex $v\in K$ the link $\lk(v)$ is combinatorially isomorphic to the boundary of a stacked polytope, then $K$ is combinatorially isomorphic to the boundary of a stacked polytope.
\end{lemma}

The following proposition seems to be new:
\begin{proposition}\label{prop:prime stress}
Let $d\geq 4$ and $K$ be the boundary complex of a prime $d$-polytope, or a prime homology $(d-1)$-sphere. Then every vertex $u\in K$ participates in a generic $d$-stress of the graph of $K$.
\end{proposition}
\begin{proof}
Let $u\in K$ be a vertex.
If $\lk(u)$ is not stacked, then $\g_2(\lk(u))>0$ by Theorem \ref{thm:stacked,Kalai}, hence $\Ker(\Rig(\lk(u)_{\leq 1},d-1))\neq \{0\}$ by 
Lemma \ref{lem:kalai rigid spheres}. So, the Cone Lemma implies that $u$ participates in a generic $d$-stress of the graph of the closed star of $u$ in $K$, hence, of $K$.
Similarly, if there exists an edge $e$ in $\antist(u)-\lk(u)$ whose two vertices are in $\lk(u)_0$, then as by Lemma \ref{lem:kalai rigid spheres} 
$\lk(u)_{\leq 1}$ is generically $(d-1)$-rigid, we get  $\Ker(\Rig(\lk(u)_{\leq 1}\cup \{e\},d-1))\neq \{0\}$. Thus, by
the Cone Lemma $u$ participates in generic $d$-stress of the graph of 
$(\{u\}*(\lk(u)\cup \{e\}))_{\leq 1}$, hence also of the graph of
$K$.

Thus, assume that (i) $\lk(u)$ is stacked and that (ii) $\antist(u)-\lk(u)$ contains no edges with both ends in $\lk(u)_0$.
Recall that the missing faces in the boundary complex of a stacked $n$-polytope different from a simplex have dimension either $1$ or $n-1$.
We now show that no facet of $\antist(u)$ has all of its vertices in $\lk(u)$; in particular $(\antist(u)-\lk(u))_0$ is nonempty. 
We show more: if $F$ is a face in $\antist(u)-\lk(u)$ then $F$ has a vertex which is not in $\lk(u)$. 
Indeed, a minimal face $F'$ in $\antist(u)-\lk(u)$ all of its vertices are in $\lk(u)$ must have size $>2$ by (ii), its boundary is contained in $\lk(u)$ by minimality, hence by (i) $F'$ is a missing facet of $\lk(u)$. Thus $F'\cup \{u\}$ is a missing facet of $K$, contradicting the fact that $K$ is prime. Thus, such $F'$ does not exist.
In particular any facet of $\antist(u)$ has a vertex not in $\lk(u)$.
We conclude that  $\antist(u)=\cup_{v\in (\antist(u)-\lk(u))_0}\overline{\st(v)}$.

By Lemma \ref{lem:kalai rigid spheres} and the Cone Lemma, for any vertex $v\in(\antist(u)-\lk(u))_0$ the closed star $\overline{\st(v)}$ is generically $d$-rigid; and it is contained in $\antist(u)$. 
Next we will show that the induced graph ($1$-skeleton) in $K$ on the vertex set 
$(\antist(u)-\lk(u))_0$ is connected. That being shown, we can totally order the vertices of $(\antist(u)-\lk(u))_0$ such that the induced graph on any initial segment is connected, say by $v_0, v_1,..., v_t$. Thus, for any $1\leq i\leq t$ the intersection $\overline{\st(v_i)}\cap (\cup_{j<i}\overline{\st(v_j)})$ contains a facet of $K$, hence by repeated application of the Gluing Lemma  
$\antist(u)=\cup_{v\in (\antist(u)-\lk(u))_0}\overline{\st(v)}$ has a generically $d$-rigid graph.

To show that the induced graph $G$ on $(\antist(u)-\lk(u))_0$ is connected, assume the contrary. 
As $\antist(u)$ is a homology ball, by Proposition \ref{prop:facet-connected} it is facet connected. Recall that any facet of $\antist(u)$ has a vertex in $(\antist(u)-\lk(u))_0$.
By assumption, there are two connected components $A$ and $B$ in $G$ and two facets $F_A$ and $F_B$ of $\antist(u)$ intersecting in a face of dimension $d-2$ such that the (one vertex in the) intersection of $F_A$ (resp. $F_B$) with  $(\antist(u)-\lk(u))_0$ is contained in $A$ (resp. $B$). Then $F=F_A\cap F_B$ belongs to $\antist(u)-\lk(u)$ and all of its vertices are in $
lk(u)$, which we showed is impossible.

Recall that $\antist(u)$ has a generically $d$-rigid graph.
Adding the edges with $u$ can increase the rank of the rigidity matrix of $\antist(u)_{\leq 1}$ by at most $d$. As $K$ is prime, $u$ has at least $d+1$ neighbors, hence the edges with $u$ contribute to the kernel of the rigidity matrix, i.e. $u$ participates in a generic $d$-stress of $K_{\leq 1}$.
\end{proof}
\begin{problem}\label{prob:edge in stress}
Do the conditions in Proposition \ref{prop:prime stress} imply that every \emph{edge} of $K$ participates in some generic stress? Equivalently, is it true that for any edge $e$ in $K$, the graph $K_{\leq 1}-\{e\}$ is generically $d$-rigid? 
\end{problem}
Note that Theorem \ref{thm:Main} implies that for a prime polytope $P$ with $\g_2(\partial P)=1$ every \emph{edge} participates in the nontrivial stress (which is unique up to nonzero scalar multiple).

\section{Proof of Theorem \ref{thm:Main}}\label{sec:Proofs}
Theorem \ref{thm:Main} is proved by induction on dimension, based on the lemmata below.
The following proposition allows the inductive step:
\begin{proposition}\label{prop:inductive step}
Let $d>4$ and $K$ be the boundary complex of a prime $d$-polytope, or a prime homology $(d-1)$-sphere, with $\g_2(K)=1$.
Then, there exists a vertex $u\in K$ such that $\lk(u)$ satisfies:

(a) $g_2(\lk(u))=1$.

(b) $\lk(u)_0=K_0-\{u\}$.

(c) $\lk(u)$ is prime.
\end{proposition}
\begin{proof}
As $g_2(K)>0$, $K$ is not stacked and by Lemma \ref{lemma:kalai stacked links} there exists a vertex $u\in K$ whose link is not stacked. By Theorem \ref{thm:stacked,Kalai}, $\g_2(\lk(u))>0$. By the Cone Lemma and Lemma \ref{lem:kalai rigid spheres}, $\g_2(\lk(u))=\dim \Ker \Rig((\overline{\st(u)})_{\leq 1},d)\leq \dim \Ker \Rig(K_{\leq 1},d)=\g_2(K)=1$. Hence $\g_2(\lk(u))=1$, proving (a).

By (a) and the Cone Lemma, there is a nontrivial $d$-stress in $K$ in which only vertices in $\overline{\st(u)}$ participate. As $\g_2(K)=1$,
no other vertex in $K$ participate in any nontrivial stress.
By Proposition \ref{prop:prime stress}, $K_0=(\overline{\st(u)})_0$, proving (b).

Assume by contradiction that $\lk(u)$ is not prime, hence inserting all of its missing facets cuts $\lk(u)$ into at least two parts, one of which is prime and the others must be simplices (e.g. it follows from Proposition \ref{prop:prime stress}, as otherwise one gets two independent generic $(d-1)$-stresses in the graph of $\lk(u)$. Alternatively use the fact that $\g_2(L\sharp Q)=\g_2(L)+\g_2(Q)$ for a connected sum $L\sharp Q$).
Denote by $M$ the prime part, and let $w$ be a vertex in $\lk(u)-M$.
By Theorem \ref{thm:stacked,Kalai} and the Cone Lemma $(\{u\}*M)_{\leq 1}$ has a generic $d$-stress, and by Proposition \ref{prop:prime stress} $w$ participates in a generic $d$-stress of $K_{\leq 1}$. Together this implies $\g_2(K)\geq 2$, a contradiction proving (c).
\end{proof}

The following two propositions establish the base of induction, namely the case $d=4$.
\begin{proposition}\label{prop:d=4,exists missing trianle}
Let $K$ be the boundary complex of a prime $4$-polytope, or a prime homology $3$-sphere, with $g_2(K)=1$ and with a missing triangle $T$. Then $K$ is combinatorially isomorphic to the join of the boundary complexes of $T$ and a polygon.
\end{proposition}
\begin{proof}
Let $T=\{a,b,c\}$. First we show that $\lk(a)_0=K_0-\{a\}$. Note that $\lk(a)_{\leq 1}$ is generically $3$-rigid by Lemma \ref{thm:rigid stress free 2-spheres}.
This graph together with the edge $\{b,c\}$ has a generic $3$-stress.
By the Cone Lemma, $G:=(\{a\}*(\lk(a)\cup \{b,c\}))_{\leq 1}$ has a generic $4$-stress. $G$ is contained in $K_{\leq 1}$.
Assume by contradiction the existence of a vertex $v\in K_0-G_0$. By Proposition \ref{prop:prime stress}, $v$ participates in a generic $4$-stress of $K_{\leq 1}$. Such a stress is independent of the former stress that we found, resulting in $g_2(K)\geq 2$, a contradiction.

By the Cone Lemma and Lemma \ref{lem:kalai rigid spheres}, the graph of $\overline{\st(a)}$ is generically $4$-rigid and $4$-stress free. This graph contains $K_0$, hence there exists exactly one edge in $\antist(a)-\lk(a)$, which must be $\{b,c\}$.

Let $C$ denote the cycle $\lk(\{b,c\},K)$, and $\Sigma C$ denote the suspension of $C$ by $b$ and $c$, namely $\Sigma C=\{b\}*C\cup \{c\}*C$. Next we show that
$\lk(a)$ contains $\Sigma C$. We have already seen that $\lk(a)_{\leq 1}\supseteq (\Sigma C)_{\leq 1}$. If a triangle $F\in \Sigma C$ is not contained in $\lk(a)$ then $F\cup \{a\}$ is missing in $K$, contradicting that $K$ is prime. Thus $\Sigma C \subseteq \lk(a)$, hence $L:=\{a\}*\Sigma C\cup_{\Sigma C}\overline{\st(\{b,c\})}=C*\partial(T)$ is the boundary of a $4$-polytope such that $L\subseteq K$.
Note that by Alexander duality a homology $d$-sphere cannot \emph{strictly} contain another homology $d$-sphere, hence $K=L=C*\partial(T)$ for $K$ a homology sphere.
\end{proof}

\begin{proposition}\label{prop:d=4,no missing trianle}
If $K$ is the boundary complex of a prime $4$-polytope, or a prime homology $3$-sphere, with $g_2(K)=1$, then $K$ has a missing triangle.
\end{proposition}
\begin{proof}
Assume by contradiction that $K$ has no missing triangles. As $K$ is prime, all of its missing faces are edges (a complex with this property is called \emph{clique complex}). Let $u\in K_0$.
Note that for every vertex $w\in \lk(u)$, $\lk(w,\lk(u))$ is a cycle of length at least $4$, as $K$ is a clique complex. Similarly, there is a vertex in $\antist(w,\lk(u))-\lk(w,\lk(u))$, as $\lk(u)$ is a non-acyclic clique complex.

Let $v\in \lk(u)$ and $I:=\antist(v,\lk(u))_0-\lk(v,\lk(u))_0$. Then $0<|I|\leq |\lk(u)_0|-5$.
Note that $K$ contains no face of the form $\{v\}\cup F$ where $F\in \antist(v,\lk(u))-\lk(v,\lk(u))$, as $K$ is a clique complex. Let $K':=(K-\st(u))\cup \{v\}*\antist(v,\lk(u))$.
Note that any edge in a clique complex satisfies the link condition. Thus,
if $K$ is a homology $3$-sphere then Proposition \ref{prop:edge-contraction-homology-sphere} asserts that $K'$ is a homology $3$-sphere as well.
(In fact, for $K$ the boundary of a simplicial polytope, $K$ and $K'$ are piecewise linearly homeomorphic \cite[Theorem 1.4]{Nevo-VK}, but we do not know whether $K'$ must be polytopal.)

Next, let us verify that $K'$ is not stacked, by showing that $K'$ is prime. As $K$ is a clique complex, any face in $K'-K$ contains an edge $\{v,i\}$ for some $i\in I\neq \emptyset$ and its vertices are contained in $\lk(u,K)_0$. Together with the fact that $K$ is prime, this implies that all the vertices of a missing tetrahedron of $K'$ must lie in $\lk(u,K)_0$ and contain $v$. However, the induced complex in $K'$ on $\lk(u,K)_0$ is a cone (over $v$), hence contains no missing tetrahedra. In particular, $K'$ is not stacked (clearly $K'$ is not the $4$-simplex, say as $K$ has at least $8$ vertices, by an easy induction on the dimension, noticing that the link of a vertex in a clique homology sphere is also a clique homology sphere).

On the other hand, $\g_2(K')=\g_2(K)-(|\lk(u,K)_0|-|I|)+4 \leq 1-5+4=0$.
This contradicts Theorem \ref{thm:stacked,Kalai}.
\end{proof}

\begin{remarks}
(1) Alternatively, Proposition \ref{prop:d=4,no missing trianle} can be proved via the Charney-Davis conjecture \cite{Charney-Davis} which was proved for homology $3$-spheres by Davis and Okun \cite{Davis-Okun}. It asserts that for a homology clique $3$-sphere $K$, $\g_2(K)-(f_0(K)-5)+1\geq 0$. In our case ($\g_2(K)=1$) we get $f_0(K)\leq 7$. As $K$ is a clique sphere, it contains at least as many vertices as the octahedral $3$-sphere, e.g. \cite[Theorem 1.1]{Meshulam-DomAndHom}, i.e. $8$ vertices; a contradiction.

(2) A homology-free though more involved proof of Proposition \ref{prop:d=4,no missing trianle} for polytopes is presented in the second author's Master thesis \cite{NovinskyMaster}.
\end{remarks}

\emph{Proof of Theorem \ref{thm:Main}}:
By Propositions \ref{prop:d=4,exists missing trianle} and \ref{prop:d=4,no missing trianle} the assertion holds for $d=4$. For $d>4$, by Proposition \ref{prop:inductive step} there exists a vertex $u\in K$ such that $K_0=\{u\}\cup \lk(u)_0$ and the conditions of Theorem \ref{thm:Main} hold for $\lk(u)$, thus by induction also the conclusion of Theorem \ref{thm:Main} holds for $\lk(u)$. Clearly, $\antist(u)-\lk(u)$ is nonempty, and any face in $\antist(u)-\lk(u)$ must contain a missing face of $\lk(u)$.
By the Cone Lemma, $g_2(\overline{\st(u)})=1=g_2(K)$, hence all the edges in $\antist(u)$ are already in $\lk(u)$. Note that the missing faces in a join are the faces which are missing in one of its components.

Case 1: $\lk(u)=\partial\sigma * C$ for a $(d-3)$-simplex $\sigma$ and a cycle $C$ of length $\geq 4$.
Since any face $F$ of $\antist(u)-\lk(u)$ contains no missing edge of $\lk(u)$, $F$ has to contain a missing face of $\partial\sigma$, namely $\sigma$.
Therefore $\sigma\in K$. As $\lk(\sigma,K)$ is a cycle and is contained in $C$, $\lk(\sigma,K)=C$. Thus, $\partial(\sigma \cup \{u\})* C \subseteq K$, and by Alexander duality $\partial(\sigma \cup \{u\})* C = K$.

Case 2: $\lk(u)=\partial\sigma * \partial\tau$ for simplices $\sigma$ and $\tau$ whose dimensions add up to $d-1$.
Then $\sigma\in K$ or $\tau\in K$ (and we shall see that exactly one of them is in $K$).
If $\sigma\in K$, then as $\lk(\sigma,K)$ is a boundary of a $(\dim \tau)$-polytope /a homology $(\dim \tau -1)$-sphere and is contained in $\partial\tau$ we must have $\lk(\sigma,K)=\partial\tau$. Then $\partial(\sigma\cup \{u\})*\partial\tau \subseteq K$ and by Alexander duality $\partial(\sigma\cup \{u\})*\partial\tau = K$.

Otherwise, $\tau\in K$ and by a similar argument $\partial(\tau\cup \{u\})*\partial\sigma = K$. $\square$

\begin{thmRemark}
Let $d\geq 4$ and $K$ be a $(d-1)$-dimensional combinatorial manifold without boundary and with $g_2(K)=1$. Then $K$ is homeomorphic to a sphere, and hence Theorem \ref{thm:Main} applies to $K$.
\end{thmRemark}
\begin{proof}
If $K$ is not prime then either it has a connected sum decomposition $K=L\sharp Q$ such that $\g_2(L)=1$ and $\g_2(Q)=0$ or it is obtained by handle forming from another combinatorial manifold without boundary $K'$ (i.e. by combinatorially identifying two disjoint closed facets of $K'$ and deleting their interior).
Here we used the fact that $d\geq 4$; see \cite{BagchiDattaSphereBundels} for a proof of this fact.

In the first case, by Theorem \ref{thm:stacked,Kalai} $Q$ is a stacked sphere and by induction on the number of vertices $L$ is homeomorphic to a sphere, and we are done.
In the second case, $\g_2(K')=\g_2(K)-\binom{d+1}{2}<0$ contradicting Theorem \ref{thm:stacked,Kalai}.

Assume that $K$ is prime. If there exists another simplicial complex $M$ which is PL-homeomorphic to $K$, and with smaller $\g_2$ value, then $M$ is a stacked sphere, hence $K$ is a PL-sphere. Otherwise, Swartz \cite{SwartzTopoFiniteness} showed that $K$ has at most $d+2$ vertices and hence $K$ is a PL-sphere \cite{BarnetteGannon}.
\end{proof}

It is natural to ask for a characterization of (prime) simplicial polytopes with a given $g_2$. First, observe the following:
\begin{observation}
Let $\g$ be the $g$-vector of a simplicial $d$-polytope with $d\geq 4$ and $\g_2>0$. Then there exists a prime $d$-polytope whose $g$-vector agrees with $\g$ except maybe in the $\g_1$ entry.
\end{observation}
\begin{proof}
The connected sum $L\sharp Q$ of two $d$-polytopes $L$ and $Q$ satisfies $\g_2(L\sharp Q)=\g_2(L)+\g_2(Q)$ and $\g_1(L\sharp Q)=\g_1(L)+\g_1(Q)+1$ (where $\g_1(L):=f_0(L)-d-1$). There exists a unique positive integer $c$ such that $\binom{c}{2}<\g_2\leq \binom{c+1}{2}$. By the sufficiency part of the $g$-theorem \cite{Billera-Lee} there exists a simplicial polytope $P$ with $\g_1(P)=c$ and $\g_i(P)=\g_i$ for any $2\leq i$. By the necessity part of the $g$-theorem \cite{Stanley:NumberFacesSimplicialPolytope-80} any simplicial polytope $P$ with $g_2(P)=\g_2$ satisfies $\g_1(P)\geq c$.
In particular, the minimality of $c$ implies that if $P=L\sharp Q$ then none of $L,Q$ is a simplex, hence w.l.o.g. $\g_1(L),\g_1(Q)>0$. The necessity part again implies $\g_2(L)\leq \binom{\g_1(L)+1}{2}$ and $\g_2(Q)\leq \binom{\g_1(Q)+1}{2}$, hence $\g_2(P)=\g_2(L)+\g_2(Q)<\binom{\g_1(L)+\g_1(Q)+1}{2}=\binom{c}{2}$, a contradiction.
\end{proof}

Next, from Example \ref{ex:big g1} it follows that:
\begin{corollary}\label{rem:big g1 with g2=b}
For every integer $d\geq 4$ and integer $b\geq 1$, there are prime $d$-polytopes $P$ with $\g_2(\partial P)=b$ and $\f_0(\partial P)$ arbitrarily large.
\end{corollary}
\begin{proof}
Note that for a prime $d$-polytope with $d\geq 4$ performing a stellar subdivision at a \emph{ridge} $F$ (i.e. $F$ has dimension $d-2$) results in a prime $d$-polytope $P'$ with $\g_2(\partial P')=\g_2(\partial P)+1$ and $\f_0(\partial P')=\f_0(\partial P)+1$. Indeed $P'$ is prime as the  missing faces in $P'$ which are not missing in $P$ are $F$, some edges (with the new vertex in $P'$) and possibly a triangle consisting of the new vertex and the two vertices in the symmetric difference between the two facets of $P$ containing $F$.

Thus, by repeating the operation of stellar subdivision over a ridge $b-1$ times, starting with the polytope in Example \ref{ex:big g1} for $n$ large, gives a polytope as claimed.   
\end{proof}

\begin{problem}
Characterize the prime polytopes with $\g_2=2$.
\end{problem}

\begin{acknowledgement}
We deeply thank Prof. Gil Kalai for many inspiring and helpful discussions. Part of the results presented here are part of the Master thesis of the second author, done under Kalai's supervision. We also thank Ed Swartz and Uli Wagner for their helpful comments, and the anonymous referees for their suggestions about the presentation.
\end{acknowledgement}
\bibliographystyle{amsplain}
\bibliography{biblio,ubt,topology}

\end{document}